\documentclass[11pt,leqno]{article}
\usepackage{amsmath, amscd, amsthm, amssymb, graphics, xypic, mathrsfs}
\usepackage{fancyhdr}
\usepackage{helvet}
\usepackage{hyperref}
\hypersetup{backref}

\newcommand{\setof}[1]{\{ #1 \}}

\newcommand{\Spec}{\operatorname{Spec}}

\newcommand{\tensor}{\otimes}

\newcommand{\Hom}{\operatorname{Hom}}
\newcommand{\isomto}{\stackrel{\sim}{\longrightarrow}}
\newcommand{\Res}{\operatorname{Res}}

\renewcommand{\O}{{\mathcal O}}

\newcommand{\F}{{\mathcal F}}

\newcommand{\Z}{{\mathbb Z}}

\newcommand{\ga}{{\mathbb G}_{\bf a}}
\newcommand{\gm}{{\mathbb G}_{\bf m}}
\newcommand{\aone}{{\mathbb A}^1}

\theoremstyle{plain}
\newtheorem{thm}{Theorem}[section]
\newtheorem{lem}[thm]{Lemma}
\newtheorem{cor}[thm]{Corollary}
\newtheorem{prop}[thm]{Proposition}

\theoremstyle{definition}
\newtheorem{defn}[thm]{Definition}

\theoremstyle{remark}
\newtheorem{rem}[thm]{Remark}
\newtheorem{ex}[thm]{Example}

\numberwithin{equation}{section}

\newcommand{\shrinkmargins}[1]{
  \addtolength{\textheight}{#1\topmargin}
  \addtolength{\textheight}{#1\topmargin}
  \addtolength{\textwidth}{#1\oddsidemargin}
  \addtolength{\textwidth}{#1\evensidemargin}
  \addtolength{\topmargin}{-#1\topmargin}
  \addtolength{\oddsidemargin}{-#1\oddsidemargin}
  \addtolength{\evensidemargin}{-#1\evensidemargin}
  }

\shrinkmargins{.6}

\begin{document}
\pagestyle{fancy} 
\renewcommand{\sectionmark}[1]{\markright{\thesection\ #1}}
\fancyhead{}
\fancyhead[LO,RE]{\bfseries\footnotesize\thepage}
\fancyhead[LE]{\bfseries\footnotesize\rightmark}
\fancyhead[RO]{\bfseries\footnotesize\rightmark}
\chead[]{}
\cfoot[]{}
\setlength{\headheight}{1cm} 

\author{Aravind Asok}

\title{{\bf {Equivariant Vector Bundles on \\ Certain Affine $G$-Varieties}}}
\maketitle

\begin{abstract}
We give a concrete description of the category of $G$-equivariant vector bundles on certain affine $G$-varieties (where $G$ is a reductive linear algebraic group) in terms of linear algebra data.  
\end{abstract}

\begin{footnotesize}
\tableofcontents
\end{footnotesize}

\section{Introduction}
For the duration of the paper we work over a fixed algebraically closed field $k$ of characteristic zero; all algebraic varieties and algebraic groups considered herein will be defined over $k$.  The word variety will be synonymous with integral, separated scheme of finite type over $k$; however, we will consider possibly disconnected algebraic groups.  Given a (possibly disconnected) linear algebraic group $G$, a variety $X$ equipped with an algebraic, left $G$-action will be called a $G$-variety.  

A $G$-variety $X$ is called {\em prehomogeneous} if $G$ acts on $X$ with a Zariski dense orbit.  In this note, we will study the problem of concretely describing the category ${\sf Vec}^G(X)$ of $G$-equivariant vector bundles on prehomogeneous $G$-varieties $X$.  The simplest instance of this problem is the case when $X$ is a homogeneous $G$-variety.  In this situation, the choice of a $k$-point $x \in X(k)$ provides an equivalence of categories $x^*: {\sf Vec}^G(X) \isomto {\sf Vec}^{G_x}(\Spec k) = {\sf Rep}(G_x)$, where $G_x$ is the stabilizer group of $x$ in $G$.  

We study the general problem as follows.  For a prehomogeneous $G$-variety $X$, let $U \subset X$ denote the dense $G$-orbit in $X$.  Restriction to $U$ defines a faithful functor ${\sf Vec}^G(X) \longrightarrow {\sf Vec}^G(U)$.  Choice of a point $x \in U(k)$ determines an equivalence of categories $x^*: {\sf Vec}^G(U) \longrightarrow {\sf Rep}(G_x)$ where $G_x$ denotes the stabilizer group in $G$ of the point $x$.  We will therefore try to describe ${\sf Vec}^{G}(X)$ in terms of representations of $G_x$ equipped with additional structure.  

The guiding example is the case $X = \aone$ and $G = \gm$.  In this case, the stabilizer group of any point in the open $\gm$-orbit is trivial.  For any $\gm$-equivariant vector bundle ${\mathcal V}$ on $X$, consider the vector space $V = \Gamma(\gm,{\mathcal V}|_{\gm})^{\gm}$, i.e. the $\gm$-invariant sections over the open subset $\gm \subset \aone$.  Choose a coordinate $x$ on $\aone$.  We define a decreasing filtration on $V$ by setting $F^p(V)$ to be the subspace of sections $v \in V$ such that $x^{-p}v$ extends to a section of ${\mathcal V}$ over $X$, i.e. by ``vanishing order of sections."  This construction defines a functor $\Phi: {\sf Vec}^\gm(\aone) \longrightarrow {\sf Filt}({\sf Vect}_k)$.  The following result is due to many authors (see Theorem \ref{lem:affineline}).   

\begin{thm}
The functor $\Phi: {\sf Vec}^{\gm}(\aone) \longrightarrow {\sf Filt}({\sf Vect}_k)$ is an equivalence of categories.  
\end{thm}

 If $G$ is a (possibly disconnected) reductive linear algebraic group, a {\em spherical} $G$-variety is a normal $G$-variety $X$ for which some Borel subgroup $B \subset G$ acts on $X$ with a Zariski dense orbit.  For a general introduction to the theory of spherical varieties, the reader may consult the article \cite{Kn3} and the notes \cite{Br}.  

Suppose now $X$ is a spherical $G$-variety.  Let $U \subset X$ denote the open $G$-orbit, let $x \in U(k)$ and let $G_x$ be the stabilizer in $G$ of $x$.  In Theorem \ref{thm:fullyfaithful} we construct a fully-faithful functor from ${\sf Vec}^G(X)$ into a category of multi-filtered representations of $G_x$ where the number of filtrations is precisely the number of codimension $1$ orbits of $G^0$ (the identity connected component of $G$).  When $X$ is an affine, spherical $G$-variety, the essential image of the aforementioned functor can be explicitly determined (see Theorem \ref{thm:affinefp}).  

\subsubsection*{Relation to Previous Work}
A concrete description of the category of equivariant vector bundles on  affine or smooth complete toric varieties was obtained by Klyachko (see \cite{Kl1} Thm 2.2.1).  Recently, a paper of Syu Kato \cite{Ka} significantly generalizes Klyachko's results in a slightly different direction.  For $G$ a connected reductive group, Kato has given a description of the category of $G \times G$-equivariant vector bundles on a class of spherical varieties consisting of certain smooth $G \times G$-equivariant partial compactifications of the group $G$. 

The natural question that arises is whether it is possible to give a concrete description of the category of $G$-equivariant vector bundles on an arbitrary (normal) prehomogeneous $G$-variety. 

\subsubsection*{Overview}
In section 2, we recall various notations related to varieties with group action and consider some preliminary structural results about equivariant vector bundles.  In section 3, we discuss in detail the case of $\gm$-equivariant vector bundles on the affine line.  In section 4, we recall some structure theory of spherical varieties.  In section 5, we state and prove the main theorem.  Finally, in section 6 we consider some examples including the space of binary quadratic forms viewed as a $GL_2$-variety and the space of $2 \times 2$-matrices over $k$ viewed as a $GL_2 \times GL_2$-variety.

\subsubsection*{A Note on Category Theoretic Terminology}
Suppose ${\mathscr C}$ and ${\mathscr D}$ are two categories and $F:{\mathscr C} \longrightarrow {\mathscr D}$ is a functor.  The {\em essential image} of $F$ is the set of all objects $d$ in ${\mathscr D}$ for which there exists an object $c \in {\mathscr C}$ and an isomorphism $d \cong F(c)$.  The functor $F$ is called {\em essentially surjective} if for every object $d \in {\mathscr D}$ there exists an object $c \in {\mathscr C}$ and an isomorphism $d \cong F(c)$.  The functor $F$ is an equivalence of categories if and only if it is fully faithful and essentially surjective (see e.g. \cite{MacLane} Ch. IV \S 4 Thm. 1); in other words $F$ admits a quasi-inverse.  

\subsubsection*{Acknowledgements}
I would like to offer my deepest thanks to Bob MacPherson for his guidance, support, energy and, most of all, for trying to teach me how to become a mathematician.  Most of the work exposed below was part of my thesis written under his direction.  Special thanks go to James Parson for great discussions and advice which clarified many points below.  I would also like to thank Mark Goresky for his advice.  Finally, an apology is in order to those who have suffered through the (many) previous versions of this paper; I sincerely thank you for your comments and what merit this paper has is due in no small part to your efforts.

\section{Notation and Preliminaries}
We first recall the definitions of various classes of varieties with group action.  For the rest of this section, $G$ will denote a (possibly disconnected) linear algebraic group.  Suppose $X$ is a $G$-variety; recall that $X$ is called {\em prehomogeneous} if $G$ acts with a Zariski dense orbit and $X$ is called {\em simple} if $G$ acts with a unique closed orbit.  Finally, $X$ is called {\em (one) fixed-pointed} if $G$ acts with a unique closed orbit that is a $G$-fixed-point.  

We denote by ${\sf Rep}(G)$ the category of (finite dimensional) $k$-rational representations of $G$, and by ${\sf Vect}_k$ the category of finite dimensional $k$-vector spaces.  For a subgroup $H \subset G$, we write $\Res^G_H: {\sf Rep}(G) \longrightarrow {\sf Rep}(H)$ for the functor which views a $G$-representation as an $H$-representation by restriction of the action.  Representations $\rho: G \longrightarrow GL(V)$ will be denoted $(V,\rho)$ or $V_\rho$ and, when no confusion will arise, we will suppress $\rho$ and refer to ``the representation $V$ of $G$."

Our main objects of study will be categories of equivariant vector bundles, let us fix notation for these.  For a $G$-variety $X$, we let ${\sf Vec}^G(X)$ denote the category of finite rank $G$-equivariant locally free sheaves  on $X$ and will refer to these objects as vector bundles.  If we consider $\Spec k$ with the trivial (left) $G$-action, the structure map $s: X \longrightarrow \Spec k$ is a $G$-equivariant map.  By definition a $G$-equivariant vector bundle on $\Spec k$ is a representation of $G$.  Given a representation $W$ of $G$, a $G$-equivariant vector bundle of the form $s^*W$ will be called a {\em constant} $G$-equivariant vector bundle.

\subsubsection*{General Remarks on Categories of Equivariant Vector Bundles}
For the remainder of this section $G$ will denote a (possibly disconnected) linear algebraic group.  Let $X$ be a $G$-variety.  Suppose $U \subset X$ is a $G$-stable open dense subvariety and let $j: U \longrightarrow X$ denote the inclusion map.  Then $j$ gives rise to the pull-back functor $j^*: {\sf Vec}^G(X) \longrightarrow {\sf Vec}^G(U)$.  Given a vector bundle ${\mathcal V}$ on $U$, we will write ${\mathcal V}|_U$ for the vector bundle $j^*{\mathcal V}$.  Since any morphism of vector bundles is determined uniquely by its restriction to the generic point, we have the following result.  

\begin{lem}
\label{lem:faithful}
The functor $j^*: {\sf Vec}^G(X) \longrightarrow {\sf Vec}^G(X)$ is faithful.  In other words, given two $G$-equivariant vector bundles ${\mathcal V}_1$ and ${\mathcal V}_2$ the map 
\begin{equation}
j^*: \Hom_{{\sf Vec}^G(X)}({\mathcal V}_1,{\mathcal V}_2) \longrightarrow \Hom_{{\sf Vec}^G(U)}({\mathcal V}_1|_{U},{\mathcal V}_2|_U)
\end{equation}
is injective.  
\end{lem}

Suppose $X$ is a normal variety, and suppose $U \subset X$ is an open dense subvariety of $X$ whose complement has codimension $\geq 2$. Given any regular function $f$ on $U$, by normality, there exists a unique regular function $f'$ on $X$ with $f'|_U = f$.  If $j: U \hookrightarrow X$ denotes the inclusion map, pull-back determines a morphism of sheaves $\O_X \longrightarrow j_*\O_U$; this map is an isomorphism.  

Locally on $U$, a morphism of locally free sheaves $\varphi: {\mathcal V}_1 \longrightarrow {\mathcal V}_2$ is given by a matrix of regular functions which extend to $X$.  More precisely, cover $X$ by open affine subvarieties $U_\alpha$ such that ${\mathcal V}_1|_{U_{\alpha}}$ and ${\mathcal V}_2|_{U_\alpha}$ are free.  A morphism $\varphi: {\mathcal V}_1|_{U \cap U_\alpha} \longrightarrow {\mathcal V}_2|_{U \cap U_\alpha}$ is given by a matrix of regular functions on $U \cap U_\alpha$ which extend to regular functions on $U_\alpha$.  Since $U$ is dense in $X$, these local extensions glue to give a global extension $\varphi:{\mathcal V}_1 \longrightarrow {\mathcal  V}_2$.  Since $G$-equivariance can be checked on a $G$-stable open dense subset, we have the following result.  

\begin{lem}
\label{lem:dense}
Suppose $X$ is a normal $G$-variety and let $U$ denote a $G$-stable open dense subvariety whose complement has codimension $\geq 2$.  Let $j: U \longrightarrow X$ denote the inclusion map.  Then
\begin{equation}
\begin{split}
j^*:{\sf Vec}^G(X) &\longrightarrow {\sf Vec}^G(U) \\
\end{split}
\end{equation}  
is a fully-faithful functor.  In particular, for any ${\mathcal V}_1,{\mathcal V}_2 \in {\sf Vec}^G(X)$,  any $G$-equivariant morphism $\varphi: {\mathcal V}_1|_U \longrightarrow {\mathcal V}_2|_U$ extends uniquely to a $G$-equivariant morphism $\varphi: {\mathcal V}_1 \longrightarrow {\mathcal V}_2$.  
\end{lem}

\subsubsection*{Some Remarks on the General Prehomogeneous Case}
\begin{ex}
\label{ex:homogeneousspace}
Suppose $X$ is a homogeneous $G$-space.  Choose a point $x \in X(k)$, and let $G_x$ denote the stabilizer group of $x$ in $G$.  This choice determines an isomorphism $X \cong G/G_x$.  The functor 
\begin{equation}
x^*: {\sf Vec}^G(X) \longrightarrow {\sf Rep}(G_x)
\end{equation}
is an equivalence of categories.  One can construct an explicit inverse functor as follows.

Consider the morphisms $\pi: G \longrightarrow G/G_x$ and $f: G \longrightarrow \Spec k$.  These are both smooth surjective morphisms.  If we equip $\Spec k$ with the trivial right $G_x$-action, $f$ is right $G_x$-equivariant as well.  Let $V$ be a $G_x$-representation; then $f^*V$ is a $G \times G_x$-equivariant vector bundle on $G$.  By faithfully flat descent (see \cite{SGA1} \S VIII Theorem 1.1) , $f^*V$ descends to a $G$-equivariant locally free sheaf on $G$.  The locally free sheaf obtained from $V$ in this manner will be denoted $G \wedge_{G_x} V$ and is called the {\em twisted or contracted product bundle}.  In fact, faithfully flat descent shows this construction is functorial.  Indeed, if $\varphi: V \longrightarrow V'$ is a morphism of $G_x$-modules, then the map $Id \times \varphi$ gives a morphism of $G \times {G_x}$-equivariant vector bundles on $G$ which descends to a morphism of $G$-equivariant sheaves:
\begin{equation}
G \wedge_{G_x} \varphi: G \wedge_{G_x} V \longrightarrow G \wedge_{G_x} V'.
\end{equation}

The functor 
\begin{equation}
G \wedge_{G_x} (\cdot): {\sf Rep}(G_x) \longrightarrow {\sf Vec}^G(X)
\end{equation}
is a quasi-inverse to the functor $x^*$.  
\end{ex}

\begin{lem}
\label{lem:prehomog}
Suppose $X$ is a prehomogeneous $G$-variety.  Suppose $U \subset X$ is the open dense orbit and let $x \in U(k)$; we denote also by $x$ the corresponding element of $X(k)$.  Let $G_x$ denote the stabilizer in $G$ of $x$.  Then 
\begin{enumerate}
\item[i)] The functor $x^*: {\sf Vec}^G(X) \longrightarrow {\sf Rep}(G_x)$ is faithful.  
\item[ii)] If $X$ is normal and $U \subset X$ has complement of codimension $\geq 2$, then $x^*$ is fully-faithful.  
\end{enumerate}
\end{lem}
 
\begin{proof}
Part (i) follows by combining Lemma \ref{lem:faithful} with the discussion of Example \ref{ex:homogeneousspace}.  Part (ii) follows by combining Lemma \ref{lem:dense} with the discussion of Example \ref{ex:homogeneousspace}.  
\end{proof}

\section{Filtered Vector Spaces and the Affine Line}
We devote this section to a discussion of $\gm$-equivariant vector bundles on the affine line $\aone$.  The construction given below is due to many authors; for example, one can refer to Gerstenhaber \cite{G} or Klyachko \cite{Kl1}.  

Recall that a {\em filtered vector space} is a pair $(V,F^{\bullet})$ consisting of a $k$-vector space $V$ and a decreasing sequence of  vector subspaces $F^{i}(V)$ indexed by $i \in \Z$.  We require that our filtrations be {\em finite} in the sense that $F^i(V) = \emptyset$ for all sufficiently large $i$ and $F^i(V) = V$ for all sufficiently small values of $i$. A morphism of filtered vector spaces $f: (V_1,F_1^{\bullet}) \longrightarrow (V_2,F_2^{\bullet})$  is a linear map $f: V_1 \longrightarrow V_2$ such that for every $i \in \Z$, we have $f(F^i_1(V_1)) \subset F^i_2(V_2)$.  We denote by ${\sf Filt}({\sf Vect}_k)$ the category whose objects are filtered vector spaces and whose morphisms are morphisms of filtered vector spaces.

We let $\Z-{\sf gr}({\sf Vect}_k)$ denote the category of $\Z$-graded vector spaces.  The ``associated graded" construction, denoted $gr$, gives a functor
\begin{equation}
gr: {\sf Filt}({\sf Vect}_k) \longrightarrow \Z-{\sf gr}({\sf Vect}_k).
\end{equation}
Consider $\aone$ with its usual $\gm$-action.  Choose a coordinate $x$ on $\aone$, and consider the $\Z$-graded ring $k[x]$ (with $x$ of degree $1$).  Let $\Z-{\sf gr}({\sf Mod}^f(k[x]))$ denote the category of finitely-generated, free, $\Z$-graded $k[x]$-modules.  The global sections functor gives an equivalence of categories
\begin{equation}
\Gamma: {\sf Vec}^{\gm}(\aone) \longrightarrow {\Z-{\sf gr}}({\sf Mod}^f(k[x])).
\end{equation}
There is an explicit inverse functor of ``sheaf associated to a graded module" which will be denoted $M \mapsto \widetilde{M}$.  

We now construct a functor $\Phi: {\sf Vec}^{\gm}(\aone) \longrightarrow {\sf Filt}({\sf Vect}_k)$.  Let $U = \aone - 0$ denote the open dense orbit.  Suppose now that ${\mathcal V}$ is a $\gm$-equivariant vector bundle on $\aone$.  Let $V$ denote the vector space of $\gm$-invariant sections of ${\mathcal V}|_U$, i.e. set
\begin{equation}
V = \Gamma(U,{\mathcal V}|_U)^{\gm}.
\end{equation}
Define a decreasing filtration of the vector space $V$ by 
\begin{equation}
F^p(V) = \setof{ v \in V |\; x^{-p}v \text{ extends to a section of } {\mathcal V} \text{ over } \aone}.
\end{equation}
The map which sends ${\mathcal V} \mapsto (V,F^{\bullet})$ is evidently functorial; moreover, it sends ${\mathcal V}$ of rank $n$ to an $n$-dimensional vector space $V$.  We denote by $\Phi: {\sf Vec}^{\gm}(\aone) \longrightarrow {\sf Filt}({\sf Vect}_k)$ the functor so defined.  By definition of the $\gm$-invariants functor, the vector space underlying $\Phi({\mathcal V})$ is equal to the fiber of ${\mathcal V}$ over the point $1 \in \gm$.  

The structure morphism $s: \aone \longrightarrow \Spec k$ is $\gm$-equivariant.  Given $W \in {\sf Rep}(\gm)$, consider the constant $\gm$-equivariant vector bundle $s^*W$.  The character group of $\gm$ is canonically isomorphic to $\Z$.  Therefore, we have a direct sum decomposition 
\begin{equation}
W \cong \oplus_{n \in \Z} W_n,
\end{equation}
where $W_n$ is the isotypic component of $W$ corresponding to the character $n$ of $\gm$.  For a constant vector bundle $\pi^*W$, one has invariant sections $V = \Gamma(U,\pi^*W|_U)^{\gm} = \oplus_{n \in \Z} x^{-n} W_n$.  
The filtration on $V$ is given by
\begin{equation}
F^p(V) = \oplus_{-n \geq p} x^{-n} W_n.
\end{equation}

\begin{thm}
\label{lem:affineline}
The functor
\begin{equation}
\Phi: {\sf Vec}^{\gm}(\aone) \longrightarrow {\sf Filt}({\sf Vect}_k)
\end{equation}
is fully faithful and essentially surjective.  
\end{thm}

\begin{proof}
The functor $\Phi$ is faithful by Lemma \ref{lem:faithful}.  

To see that $\Phi$ is full, we proceed as follows.  Given a $\gm$-equivariant vector bundle ${\mathcal V}$ on $\aone$, note that $\Gamma(\aone,{\mathcal V}) \subset \Gamma(U,{\mathcal V}|_U)$.  More precisely, $\Gamma(\aone,{\mathcal V})$ can be identified with the $\Gamma(\aone,\O_{\aone})$-submodule of $\Gamma(U,{\mathcal V}|_U)$ generated by $x^{-p}F^p(V)$ for all $p$.  

Fix a pair ${\mathcal V}_1,{\mathcal V}_2$ of $\gm$-equivariant vector bundles on $\aone$ with $\Phi({\mathcal V_i}) = V_i$.   If for a locally free sheaf morphism $f: {\mathcal V}_1|_U \longrightarrow {\mathcal V}_2|_U$, the map $\Phi(f): V_1 \longrightarrow V_2$ preserves the filtrations, then $f$ extends to a vector bundle morphism on $\aone$.  By density of $U$, we can check $\gm$-equivariance over $U$, and so $f$ actually extends to a $\gm$-equivariant vector bundle morphism over $\aone$.

To see that $\Phi$ is essentially surjective, it suffices to consider the constant $\gm$-equivariant vector bundles on $\aone$.  The computation just preceding the theorem shows that the functor $\Phi$ applied to constant $\gm$-equivariant vector bundles exhausts all filtered vector spaces up to isomorphism.  
\end{proof}

\begin{rem}[Rees Construction]
Implicit in our description above was an explicit quasi-inverse to the functor $\Phi$; this functor is sometimes known as the Rees construction (this construction is treated in much greater depth in \cite{G}).  Given a filtered vector space $(V,F^{\bullet})$.  We construct a $\Z$-graded $k[x]$-module by defining
\begin{equation}
\widehat{(V,F^{\bullet})} = \bigoplus_{n \in \Z} x^{-n} F^n(V).
\end{equation}
Note that the module on the right hand side is a finitely-generated, free, $\Z$-graded $k[x]$-module.  Then, the functor 
\begin{equation}
\widetilde{\widehat{(\cdot)}}: {\sf Filt}({\sf Vect}_k) \longrightarrow {\sf Vec}^{\gm}(\aone)
\end{equation}
is quasi-inverse to $\Phi$.  
\end{rem}

\begin{rem}[Constant Bundles and Gradings]
The Rees construction can be viewed as producing a flat deformation of a filtered vector space to its associated graded vector space.  To see this, let $0: \Spec k \longrightarrow \aone$ denote the $\gm$-fixed point in $\aone$.  Then we have the functor $0^*: {\sf Vec}^{\gm}(\aone) \longrightarrow {\sf Rep}(\gm)$.  Consider also the functor ${\sf r}: {\sf Rep}(\gm) \longrightarrow \Z-{\sf gr}({\sf Vect}_k)$ which reverses the grading; ${\sf r}$ is defined on objects $W \in {\sf Rep}(\gm)$ by
\begin{equation}
W = \oplus_{n \in \Z} W_n \longrightarrow \oplus_{n \in \Z} W^n
\end{equation}
where $W^n = W_{-n}$.  Consider the diagram
\begin{equation}
\xymatrix{
{\sf Vec}^{\gm}(\aone) \ar[r]^{\Phi}\ar[d]^{0^*} & {\sf Filt}({\sf Vect}_k) \ar[d]^{gr} \\
{\sf Rep}(G) \ar[r]^{{\sf r}} & \Z-{\sf gr}({\sf Vect}_k)
}.
\end{equation}
We claim there is a natural isomorphism of functors $gr \circ \Phi \stackrel{\sim}{\Longrightarrow} {\sf r} \circ 0^*$.

To see this, consider a $\gm$-equivariant vector bundle ${\mathcal V}$ on $\aone$.  Let $(V,F^{\bullet}) = \Phi({\mathcal V})$.  Map $F^p(V)$ to $0^*({\mathcal V})$ by sending $v \in F^p(V)$ to $x^{-p}v |_{0^*{\mathcal V}}$.  By definition, the kernel of the map so defined is isomorphic to $F^{p+1}(V)$ and so we obtain an injection $F^p(V)/F^{p+1}(V) \hookrightarrow 0^*{\mathcal V}$.  

Note furthermore that since $\gm$ acts by $t \mapsto t^{-p}$ on $x^{-p}v$ ($v$ is $\gm$-invariant) the image of $F^p(V)/F^{p+1}(V)$ is contained in $W_{-p} = W^p$.  Therefore, the map 
\begin{equation}
\label{eqn:graded}
gr_{F^{\bullet}}(V) \longrightarrow 0^*{\mathcal V}
\end{equation}
is an injection (since the $F^p(V)/F^{p+1}(V)$ map to independent isotypic subspaces), and the map in Equation \ref{eqn:graded} is an isomorphism as both sides are vector spaces of the same dimension.  
\end{rem}

\section{Some Structure of Spherical Varieties}
Let $G$ be a reductive group.  Recall that a normal $G$-variety $X$ is called {\em spherical} if some Borel subgroup $B \subset G$ acts on $X$ with a dense orbit.  Until we mention otherwise, we assume $G$ is {\em connected} and reductive.  

Let $X$ be a spherical $G$-variety.  Let $U \subset X$ denote the dense $G$-orbit.  Choose a point $x \in U(k)$ and let $H = G_x$ denote the stabilizer group.  For such an $X$, we now construct certain homomorphisms $\lambda: \gm \longrightarrow G$ such that $\lim_{t \longrightarrow 0} \lambda(t) \cdot x$ exists and lies in codimension $1$, $G$-orbit.  

\subsubsection*{Geometry of Elementary Embeddings}
General references for the material in this section are \cite{BP} \S 2 or \cite{BLV}.  Let $X$ be a spherical variety, let $U \subset X$ denote the dense $G$-orbit, choose $x \in U(k)$ and let $H = G_x$.  The pair $(X,x)$ will be called an {\em elementary (spherical) embedding} if $X$ is regular, and $X \backslash U$ consists of a single closed, codimension $1$, $G$-orbit.  We denote this closed orbit by $X'$.  Note that an elementary spherical embedding is simple by definition.  

Given an elementary embedding $(X,x)$, fix a Borel subgroup $B \subset G$ such that $BH$ is open and Zariski dense in $G$ (such a Borel subgroup exists by definition of sphericity).  Let $P = \setof{g \in G | gBH = BH}$.  Then $P$ is a parabolic subgroup of $G$ containing $B$ and $P$ has a dense orbit in $X$, namely the orbit $P \cdot x$.  In addition, one can show that $P$ has a dense orbit $Y'$ in $X'$; let $Y = P \cdot x \cup Y'$.  Furthermore, let $R_u(P)$ denote the unipotent radical of $P$.  

\begin{prop}[Brion-Luna-Vust \cite{BLV} \S 4.2] 
\label{prop:adaptedlevi}
Notation as in the previous paragraph, there exists a Levi-subgroup $L \subset P$ such that for any elementary spherical embedding $(X,x)$ of $G/H$ we have:
\begin{enumerate}
\item[i)] $P \cap H = L \cap H$,
\item[ii)] $P \cap H$ contains $[L,L]$,
\item[iii)] if $C$ denotes the identity connected component of $L$,  then the action of $R_u(P)$ on $Y$ induces a $P$-equivariant isomorphism $R_u(P) \times (\overline{C\cdot x} \cap Y) \isomto Y$.  
\end{enumerate}
\end{prop}

A Levi subgroup $L \subset P$ satisfying the conclusions of Proposition \ref{prop:adaptedlevi} is said to be {\em adapted} to $H$.  The orbit closure $\overline{C \cdot x}$ is normal, and since $C$ is a torus, is a toric variety for a quotient torus of $C$.

For any linear algebraic group $H$, let ${\rm X}_*(H)$ and ${\rm X}^*(H)$ denote the lattices of cocharacters and characters respectively.  For a torus $T$, we denote by $\langle \cdot , \cdot \rangle$ the natural pairing ${\rm X}_*(T) \times {\rm X}^*(T) \longrightarrow \Z$.  

Given a one parameter subgroup $\lambda \in {\rm X}_*(G)$, one can associate to $\lambda$ a valuation of $k(G/H)$ as follows.  If $f \in k[G]$ is a regular function, then there exist a finite number of $f_n \in k[G]$ such that $f(g\lambda(t)) = \sum_{n \in \Z} f_n(g) t^n$, for any $g \in G$ and $t \in \gm$.  Set $\overline{\nu}_\lambda(f) = \min \setof{n | f_n \neq 0}$.  Then $\overline{\nu}_\lambda$ extends to a valuation of $k(G)$ which is left $G$-invariant (this means $\overline{\nu}_\lambda(g\cdot f) = \overline{\nu}_\lambda(f)$ for all $g \in G$) and we write $\nu_\lambda$ for the composite of the natural inclusion $k(G/H) \hookrightarrow k(G)$ with $\overline{\nu}_{\lambda}$; this is again a left $G$-invariant valuation.  

Given an elementary embedding $(X,x)$, the codimension $1$ orbit defines a valuation of $k(G/H)$ by vanishing order along $X'$.  For an elementary spherical embedding, these two constructions of valuations can be related as follows.  

\begin{prop}[Brion-Luna-Vust \cite{BLV} \S 4.2]
\label{prop:onepsg}
Let $(X,x)$ be an elementary spherical embedding.  Let $L$ be a Levi subgroup of $P$ adapted to $H$.  Then there exists a one parameter subgroup $\mu: \gm \longrightarrow C$, unique up to multiplication by an element of ${\rm X}_*(C \cap H)$, such that $\lim_{t \rightarrow 0} \mu(t)\cdot x$ exists and lies in the open $P$-orbit in $X'$.  Moreover, the valuation $\nu_{\mu}$ is equivalent to the valuation of $k(G/H)$ defined by $(X,x)$. 
\end{prop}

\subsubsection*{Probing the Boundary of a Spherical Variety}
Suppose $X$ is a spherical $G$-variety.  Let $U \subset X$ denote the dense $G$-orbit and choose $x \in U(k)$; again, set $G_x = H$.  To each codimension $1$ $G$-orbit in $X$, we can associate an elementary embedding of $G/H$ as follows.  Let $D_i$ denote a codimension $1$ $G$-orbit and consider the spherical $G$-variety $X_i = U \cup D_i$.  Since spherical varieties are by definition normal, and the singular locus of a normal $G$-variety is a $G$-stable subset of codimension $\geq 2$, $X_i$ is a regular variety.  Therefore, $X_i$ is an elementary embedding of $G/H$.  

Fix an adapted Levi subgroup $L$ of $X_i$ (which exists by Proposition \ref{prop:adaptedlevi}).  By Proposition \ref{prop:onepsg}, for each $i$ we can find a one parameter subgroup $\mu_i: \gm \longrightarrow C$ such that $\lim_{t \rightarrow 0} \mu_i(t) \cdot x$ exists and lies in $D_i$.  In other words, $\mu_i$ gives rise to a $\gm$-equivariant morphism $\overline{\mu_i}: \aone \longrightarrow X_i$.  

\subsubsection*{Possibly Disconnected Reductive Groups}
Suppose now that $G$ is a possibly disconnected reductive group.  We let $G^0$ denote the identity connected component of $G$.  If $X$ is a spherical $G$-variety, then $X$ is also a spherical $G^0$-variety.  This means to each codimension $1$, $G^0$-stable boundary component, we can attach a one-parameter subgroup via Proposition \ref{prop:onepsg}.  Summarizing the discussion of the previous paragraphs, we are led to the following definition.  

\begin{defn}
\label{defn:boundaryonepsg}
Suppose $G$ is a (possibly disconnected) reductive group.  Suppose $X$ is a spherical $G$-variety, let $U \subset X$ be the dense $G$-orbit, fix $x \in U(k)$ and let $H = G_x$.  Let $I$ denote the (possibly empty, finite) set of codimension $1$, $G^0$-orbits $D_i$ in $X$, and let $X_i = U \cup D_i$ denote the associated elementary $G^0$-embeddings described above.  A collection of one-parameter subgroups $\mu_i: \gm \longrightarrow G$ such that $\mu_i$ satisfies the conclusion of Proposition \ref{prop:onepsg} for each $X_i$ will be called a family of {\em boundary one parameter subgroups} of $X$.  The $\gm$-equivariant maps defined by $t \mapsto \mu_i(t)\cdot x$ extend to $\gm$-equivariant morphisms denoted $\overline{\mu_i}: \aone \longrightarrow X$.  
\end{defn}

\begin{rem}
For general normal prehomogeneous $G$-varieties, one parameter subgroups having the properties stated in Proposition \ref{prop:onepsg} need not exist.  One can define the notion of an elementary embedding $(X,x)$ of a non-spherical homogeneous space.  For an elementary embedding of a non-spherical homogeneous space it need not be the case that there exists a one parameter subgroup $\mu: \gm \longrightarrow G$ such that $\lim_{t \rightarrow 0} \mu(t) \cdot x$ exists in the codimension $1$ $G$-orbit in $X$.  Indeed, the limit point of a one-parameter subgroup as in Proposition \ref{prop:onepsg} is necessarily fixed by $\gm$.  Consider $G = GL_2$ acting on the vector space $Sym^3 V$ where $V$ is the standard $2$-dimensional representation of $G$.  There is a single codimension $1$ $G$-orbit in $X$ and the stabilizer of a point in this orbit is isomorphic to an extension of the additive group $\ga$ by a finite group; this group does not contain a subgroup isomorphic to $\gm$.  
\end{rem}

\section{Equivariant Vector Bundles on Spherical Varieties}
Suppose $X$ is a spherical $G$-variety.  As in the previous section, let $U \subset X$ denote the dense $G$-orbit, fix $x \in U(k)$, let $H = G_x$.  Furthermore, choose $B,P$ (as explained before Proposition \ref{prop:adaptedlevi}), and a Levi subgroup $L$ adapted to $H$.  As in Proposition \ref{prop:adaptedlevi} let $C$ denote the identity connected component of the center of $L$.  Choose a collection $\setof{\mu_i}_{i \in I}$ of boundary one-parameter subgroups as in Definition \ref{defn:boundaryonepsg}.

\begin{defn}
Let $I$ be a (possibly empty) finite set and let $H$ be a linear algebraic group.  Define a category ${\sf Filt}(I,H)$ as follows.
\begin{itemize}
\item[$\bullet$] Objects are collections $(V,\rho,F_i^{\bullet})$ where $V$ is a $k$-vector space, $\rho: H \longrightarrow GL(V)$ is a representation of $H$ on $V$, and $F_i^{\bullet}$ is a family of filtrations of the vector space $V$ indexed by $i \in I$.  (By convention, if $I$ is empty, objects are just pairs $(V,\rho)$.)    
\item[$\bullet$] Given collections $(V_1,\rho_1,F^{\bullet}_{i,1})$, $(V_2,\rho_2,F^{\bullet}_{i,2})$, morphisms are $k$-linear maps $f: V_1 \longrightarrow V_2$ which simultaneously are $H$-module homomorphisms and preserve the filtrations. 
\end{itemize}
\end{defn}

\subsubsection*{A Fully-Faithful Functor}
Each boundary $1$-parameter subgroup $\mu_i: \gm \longrightarrow G$ defines a functor $\overline{\mu_i}^*: {\sf Vec}^{G}(X) \longrightarrow {\sf Vec}^{\gm}(\aone)$.  Composing with the functor $\Phi$ of Theorem \ref{lem:affineline}, we obtain functors $\Phi \circ \overline{\mu_i}^*: {\sf Vec}^G{X} \longrightarrow {\sf Filt}({\sf Vect}_k)$.  We also have the functor $x^*: {\sf Vec}^G(X) \longrightarrow {\sf Rep}(H)$ defined by our chosen point $x$ in the dense $G$-orbit.  Now, the map $t \mapsto \mu_i(t) \cdot x$ sends $1 \in {\gm}$ to $x$.  Therefore, we can consider the functor
\begin{equation}
\Psi: {\sf Vec}^G(X) \longrightarrow {\sf Filt}(I,H) \;\;\; {\mathcal V} \mapsto (x^*{\mathcal V},\Phi(\overline{\mu_i}^*{\mathcal V})).
\end{equation}
The first main result is the following.

\begin{thm}
\label{thm:fullyfaithful}
The functor $\Psi: {\sf Vec}^G(X) \longrightarrow {\sf Filt}(I,H)$ is fully faithful. 
\end{thm}

\begin{proof}
That $\Psi$ is faithful follows from Lemma \ref{lem:prehomog} (i).  

To show $\Psi$ is full we proceed as follows.  Given two $G$-equivariant vector bundles ${\mathcal V}_1$ and ${\mathcal V}_2$ on $X$, we show that any morphism $\varphi: \Psi({\mathcal V}_1) \longrightarrow \Psi({\mathcal V}_2)$ extends to a bundle morphism $\overline{\varphi}: {\mathcal V}_1 \longrightarrow {\mathcal V}_2$.  Note that $\varphi$ gives rise to morphism $\overline{\varphi}: {\mathcal V}_1|_U \longrightarrow {\mathcal V}_2|_U$.  

If $I$ is empty, the fact that $\Psi$ is full follows from $\ref{lem:prehomog} (ii)$.  In general, the union of $U$ with the $G$-orbits of codimension $1$ is $G$-stable and has complement of codimension $\geq 2$.  Therefore if we show that $\overline{\varphi}$ extends across $D_i$ for each $i$, it follows that it extends to a vector bundle morphism on $X$ by Lemma \ref{lem:dense}.  Let us now show that $\overline{\varphi}$ extends across $D_i$ for each $i$.  

The valuation $\nu_{\mu_i}$ is equivalent to the valuation of $X_i$ by Proposition \ref{prop:onepsg}.  In other words, given a rational function $f$ on $X_i$, $f$ is regular on $X_i$ if and only if $\overline{\mu_i}^*f$ is regular at $0$.  The morphism of vector bundles $f$ is determined by a matrix of rational functions on $X$ and we need to check that each of these rational functions is regular along $D_i$.  

Since $\varphi$ preserves the filtration $F_i^{\bullet}$, we know that $\overline{\mu_i}^*\varphi: \overline{\mu_i}^* {\mathcal V}_1 \longrightarrow \overline{\mu_i}^* {\mathcal V}_2$ extends to a morphism of $\gm$-equivariant vector bundles on $\aone$.  This means that the rational functions defining $\overline{\varphi}$ must be regular along $D_i$ since their restrictions to $\aone$ are regular at $0$.  In other words, $\overline{\varphi}$ extends to a morphism of bundles over $X_i$.  By assumption, this is true for every $i$ and hence $\overline{\varphi}$ extends across $D_i$ for every $i$.   
\end{proof}

Our goal is to determine the essential image of $\Phi$.  Analogous to the case of the affine line (see Theorem \ref{lem:affineline}) the first natural class of objects to consider are constant $G$-equivariant vector bundles.  However, in order to do this, we will need to add additional constraints to the spherical varieties under consideration.

\subsubsection*{Geometry of Affine Spherical Varieties}
Suppose $X$ is an affine spherical $G$-variety.  Since $k[X]^G = k$, $X$ must be simple by \cite{GIT} Ch. 1 Cor 1.2, since distinct closed orbits can be separated by invariant functions.  Let $Y \subset X$ denote the unique closed orbit, fix $y \in Y(k)$ and let $G_y$ be the stabilizer group in $G$ of $y$.  Note that $Y \subset X$ is affine as a closed subvariety of an affine variety is affine.  Hence $Y \cong G/G_y$ is an affine homogeneous $G$-variety.  Matsushima's theorem (see e.g. \cite{H} Thm 3.3) states that a homogeneous space $G/G_y$ is affine if and only if $G_y$ is a (possibly disconnected) reductive subgroup of $G$.  

Luna's slice theorem constrains the geometry of affine spherical $G$-varieties $X$ as follows.  By \cite{Luslice} \S 3.1 Cor. 2, there exists a closed $G_y$-stable, fixed-pointed, subvariety $Z \subset X$ such that the (multiplication) map
\begin{equation}
m: G \wedge_{G_y} Z \longrightarrow X
\end{equation}
is a $G$-equivariant isomorphism.  If $X$ is in addition smooth, then $Z$ is $G_y$-equivariantly isomorphic to a vector space with linear $G_y$-action, i.e. a $G_y$-representation.  Let $i: Z \hookrightarrow X$ denote the inclusion map.  

\begin{lem}
\label{lem:sliceequivalence}
The functor $i^*: {\sf Vec}^G(X) \longrightarrow {\sf Vec}^{G_y}(Z)$ is an equivalence of categories.  
\end{lem}

\begin{proof}
Let $i': Z \hookrightarrow G \times Z$ be the zero section of the trivial left $G$-torsor  $p_2: G \times Z \longrightarrow Z$; $p_2$ is a faithfully flat, quasi-compact morphism.  Let $q: G \times Z \longrightarrow G \wedge_{G_y} Z$, be the quotient morphism which makes $G \times Z$ into a right $G_y$-torsor over $G \wedge_{G_y} Z$; $q$ is also a faithfully flat, quasi-compact morphism.  

Let $\pi_i: (G \times Z) \times_{p_2,Z,p_2} (G \times Z) \cong G \times (G \times Z) \longrightarrow (G \times Z)$ ($i = 1,2$) denote the two canonical projections onto the first and second factors of the fiber product.  Consider the isomorphism $(G \times Z) \times_{p_2,Z,p_2} (G \times Z) \cong G \times (G \times Z)$ defined by the map $(g,z),(g',z) \mapsto (gg', g,z)$.  Under this isomorphism, the projection maps $\pi_i$ are sent to maps $G \times (G \times Z) \longrightarrow G \times Z$: the action map $(g_1, (g_2,z)) \mapsto (g_1g_2,z)$ and the projection map $(g_1,(g_2,z)) \mapsto (g_2,z)$.  Therefore, for a vector bundle $\F$ on $G \times Z$, specifying a descent datum for $p_2$ is exactly equivalent to specifying a $G$-equivariant structure on $\F$.  Similarly, for a vector bundle $\F$ on $G \times Z$, giving a descent datum for the morphism $q$, is exactly equivalent to giving a $G_y$-equivariant structure on a bundle over $\F$.    

By faithfully flat descent (see \cite{SGA1} \S VIII Theoreme 1.1), we see that $i'^*: {\sf Vec}^{G \times G_y}(G \times Z) \longrightarrow {\sf Vec}^{G_y}(Z)$ is an equivalence of categories and that $q^*: {\sf Vec}^{G}(G \wedge_{G_y} Z) \longrightarrow {\sf Vec}^{G \times G_y}(G \times Z)$ is an equivalence of categories as well.  In each case, the equivariant structure gives, by definition, a descent datum.  Now, we just note that the morphism $i$ factors as $qi'$.  
\end{proof}

Because of this result, we will only give a concrete description of the category of equivariant vector bundles on affine fixed-pointed spherical varieties.  

\subsubsection*{Equivariant Nakayama Lemma}
Isomorphism classes of $G$-equivariant vector bundles on fixed-pointed, spherical, $G$-varieties are particularly simple.  Indeed, Bass and Haboush deduce the following fact from their equivariant Nakayama lemma (see \cite{BaHa} Proposition 6.1).  

\begin{lem}[\cite{BaHa} \S6.5]
\label{lem:fpvars}
Suppose $X$ is a fixed-pointed, affine, spherical $G$-variety.  Let $x_0 \in X$ denote the $G$-fixed point.  Suppose ${\mathcal V}$ and ${\mathcal V}'$ are a pair of $G$-equivariant vector bundles on $X$.  Then
\begin{itemize}
\item[i)] Any $G$-equivariant morphism $\varphi_0: x_0^*{\mathcal V} \longrightarrow x_0^*{\mathcal V}'$ can be lifted to a $G$-equivariant morphism $\varphi: {\mathcal V} \longrightarrow {\mathcal V}'$.
\item[ii)] Furthermore, if $\varphi_0$ is surjective, then $\varphi$ is surjective, and
\item[iii)] If $\varphi_0$ is an isomorphism, then $\varphi$ is an isomorphism. 
\end{itemize}
\end{lem}

\begin{cor}
\label{cor:affinesphericalvars}
Let $X$ be an affine, spherical $G$-variety and let $y$ be a point in the unique closed $G$-orbit in $X$ with stabilizer group $G_y$.  Let $Z \subset X$ be a $G_y$-stable fixed-pointed variety such that $X \cong G \wedge_{G_y} Z$.  Let $s: Z \longrightarrow \Spec k$ denote the structure morphism.  Every $G$-equivariant vector bundle on $X$ is isomorphic to one of the form $G \wedge_{G_y} s^*W$ where $W$ is a finite dimensional representation of $G_y$.    
\end{cor}

\begin{proof}
We combine the equivalence of categories discussed in Lemma \ref{lem:sliceequivalence} with Lemma \ref{lem:fpvars}.  
\end{proof}

\subsubsection*{Constant Vector Bundles}
Suppose $\mu: \gm \longrightarrow G$ is a one-parameter subgroup.  Suppose $(\rho,V)$ is a representation of $G$.  The homomorphism $\mu$ factors through a torus $T \subset G$; fix such a torus.  Consider the decomposition of $(\rho,V)$ as a $T$-representation; we let $\chi$ denote a character of $T$.  We then associate a filtration $F_{\mu}^{\bullet}$ of the vector space $V$ to $\mu$ by setting
\begin{equation}
\label{eqn:constantfiltration}
F_{\mu}^i(V) =  \bigoplus_{- \langle \mu, \chi \rangle \geq i} V_\chi.
\end{equation}
For any family of $1$-parameter subgroups $\mu_i$ indexed by some finite set $I$, define a functor ${\sf C}: {\sf Rep}(G) \longrightarrow {\sf Filt}(I,H)$ by 
\begin{equation}
(\rho,V) \mapsto (V,{\sf Res}^G_H(\rho),F_{\mu_i}^{\bullet}).
\end{equation}

\begin{defn}
\label{defn:constantbundles}
Let ${\sf Filt}(I,H,G)$ denote the full subcategory of ${\sf Filt}(I,H)$ consisting of objects in the essential image of ${\sf C}$.
\end{defn}

\subsubsection*{Vector Bundles on Affine Spherical Varieties}
\begin{thm}
\label{thm:affinefp}
Let $X$ be an affine, fixed-pointed, spherical $G$-variety.  Let $\setof{\mu_i}$ be a family of boundary one-parameter subgroups of $X$ as in Definition \ref{defn:boundaryonepsg}.  Then the functor $\Psi$ of Theorem \ref{thm:fullyfaithful} gives a fully-faithful essentially surjective functor
\begin{equation}
{\sf Vec}^G(X) \longrightarrow {\sf Filt}(I,G,H)  
\end{equation}
(see \ref{defn:constantbundles}).
\end{thm}

\begin{proof}
We have to check that the filtration $F_i^{\bullet}$ on a constant bundle associated to a representation $W$ is the same as the filtration $F_{\mu_i}$ on the representation $W$.  To do this, we note the morphism $\overline{\mu_i}: \aone \longrightarrow X$ fixes its limit point.  The fibers of a constant $G$-equivariant vector bundle over any two $k$-points of $X$ are canonically identified.  Hence we have a canonical isomorphism between the fiber over the unique $G$-fixed point and the fiber over the $\gm$-fixed point; we use this isomorphism to identify the filtrations.  It suffices to note that, by Lemma \ref{lem:fpvars}, every $G$-equivariant vector bundle on $X$ is isomorphic to a constant bundle.
\end{proof}

Recall that a {\em{skeleton}} of a category ${\mathcal C}$ is a full sub-category $sk({\mathcal C})$ of ${\mathcal C}$ such that the inclusion functor is essentially surjective and in which no two distinct objects are isomorphic (see \cite{MacLane} Ch. IV \S4, p. 91).  Every category has a skeleton and the inclusion functor is an equivalence of categories.  For an affine, fixed-pointed, spherical $G$-variety, it follows from Lemma \ref{lem:fpvars} that the category whose objects are constant vector bundles on $X$ and whose morphisms are $G$-equivariant vector bundle morphisms is a skeleton for ${\sf Vec}^G(X)$.  

\begin{defn}
Fix notation as in Theorem \ref{thm:affinefp}.  Let ${\sf Filt}(I,X)$ denote the following category.
\begin{itemize}
\item[$\bullet$] Objects are $G$-representations $(\rho,W)$. 
\item[$\bullet$] Morphisms $f: (\rho_1,V_1) \longrightarrow (\rho_2,V_2)$ are $k$-linear maps $f: V_1 \longrightarrow V_2$ which are $H$-module homomorphisms and which preserve the filtrations $F_{\mu_i}$ defined above Definition \ref{defn:constantbundles}.   
\end{itemize}
\end{defn}

There is a functor $\Omega: {\sf Filt}(I,X) \longrightarrow {\sf Vec}^G(X)$ defined as follows.  An object $W \in {\sf Filt}(I,X)$ is sent to the constant $G$-equivariant vector bundle.  A morphism $f \in {\sf Filt}(I,X)$  defines a morphism of constant $G$-equivariant vector bundles by the proof of Theorem \ref{thm:fullyfaithful}; in other words, we define a morphism on the open dense $G$-orbit and extend it to $X$ by normality.    

\begin{thm}
\label{thm:rephrase}
The functor $\Omega: {\sf Filt}(I,X) \longrightarrow {\sf Vec}^G(X)$ is an equivalence of categories.  
\end{thm}

\begin{proof}
The functor $\Omega$ identifies ${\sf Filt}(I,X)$ with a skeleton of the category ${\sf Vec}^G(X)$ consisting of constant $G$-equivariant vector bundles.  
\end{proof}

\section{Examples}
For an affine $G$-variety $X$, the space of homomorphisms in the category ${\sf Vec}^G(X)$ can be related to the $G$-module structure of $k[X]$.  In what follows, we will consider some infinite rank $G$-modules which are direct limits of finite dimensional $G$-modules.  Formally, we can consider the category ${\sf Mod}(G)$ of inductive limits of finite dimensional representations.  We let $\Hom_G$ denote the space of morphisms in this category.  The following lemma will be used in connection with the main theorem to give explicit formulae for multiplicities of $G$-modules in $k[X]$. 

\begin{lem}
\label{lem:mult}
Let $X$ be an affine $G$-variety.  Suppose $V_i$ ($i = 1,2$) are finite dimensional $G$-modules.  Then there is a canonical isomorphism
\begin{equation}
\Hom_{{\sf Vec}^G(X)}(X \times V_1,X \times V_2) \cong \Hom_{G}(V_1,V_2 \tensor k[X]).
\end{equation} 
\end{lem}

\begin{proof}
By duality, specifying a $G$-equivariant bundle map $\varphi: X \times V_1 \longrightarrow X \times V_2$ is equivalent to specifying a regular map $\varphi': X \longrightarrow V_1^{\vee} \times V_2$ which is $G$-invariant.  Such a map is by definition an element of the $k$-vector space $(V_{1}^{\vee} \tensor V_2 \tensor k[X])^{G}$.  Specifying an element of this last space is in turn equivalent to specifying an element of $Hom_{G}(k,V_{1}^{\vee} \tensor V_2 \tensor k[X])$.  Using duality again, this is isomorphic to the group $Hom_{G}(V_1,V_2 \tensor k[X])$.  
\end{proof}

\begin{ex}[Binary Quadratic Forms]
\label{ex:binquadform}
Consider the $GL_2$-action on the space of binary quadratic forms: i.e. we look at the action $GL_2$ on the space of two-variable polynomials of the form $q = ax^2 + bxy + cy^2$ by substitution.  In other words, if $V$ denotes the standard $2$-dimensional representation of $GL_2$, we consider the $GL_2$-action on $Sym^2 V^\vee$.  By explicit computation, one can see that this is a spherical $G$-variety.  We completely describe the $GL_2$-module structure of $k[Sym^2 V^{\vee}]$.  

Here, $GL_2$ acts with three orbits.  If we define the discriminant of a quadratic form $q$ (as above) to be $b^2 - 4ac$, then the open orbit can be identified with non-vanishing locus of the discriminant.  The closed subvariety where the discriminant vanishes consists of two $GL_2$-orbits: the $GL_2$-fixed point corresponding to the quadratic form where $a = b = c = 0$ and its complement.  

The quadratic form $xy + y^2$ lies in the open $GL_2$-orbit as it has non-zero discriminant.  Define 
\begin{equation}
\mu(t) = \begin{pmatrix}
t & 0 \\
0 & 1
\end{pmatrix}.
\end{equation}
Note that $\lim_{t \rightarrow 0} \mu(t) \cdot (xy + y^2)$ exists in $Sym^2 V^{\vee}$ and is equal to $y^2$.  This one parameter subgroup vanishes to order $1$ on the discriminant locus and satisfies the conclusion of Proposition $\ref{prop:onepsg}$.  

The stabilizer of the quadratic form $xy + y^2$ can be identified with the subgroup $H$ of $GL_2$ whose elements are of the form
\begin{equation}
\begin{pmatrix}
t & t^{-1} - t \\
0 & t^{-1}
\end{pmatrix}
\end{equation}
for a non-zero $t$.  Conjugation by the matrix $\begin{pmatrix} 1 & 1 \\ 0 & 1\end{pmatrix}$ in $GL_2$ gives an isomorphism $H \cong \gm$; here $\gm$ is identified with a maximal torus in $SL_2$.  

Every irreducible representation of $GL_2$ is isomorphic to one of the form $Sym^{n} V \tensor (\wedge^2 (V))^{\tensor m}$ for some positive integer $n$ and an arbitrary integer $m$.  Under the maximal torus of diagonal matrices $diag(t_1,t_2)$ in $GL_2$, this representation is $n+1$-dimensional and has weights 
\begin{equation}
\label{eqn:weights}
(t_1^{n+m}t_2^{m},t_1^{n-1+m}t_2^{m+1},\ldots,t_1^{m}t_2^{n-1+m},t_1^{m}t_2^{n+m}).
\end{equation}
We refer to this representation as ``the representation $(n,m)$."  

For the representation $(n,m)$, $\mu(t)$ acts with weights $t^{m+n},t^{m+n-1},\ldots,t^{m+1},t^{m}$.   Hence the filtration on the vector space $W$ underlying such a representation is given by (see Equation \ref{eqn:constantfiltration})
\begin{equation}
W_{\geq -(m+n)} \supset W_{\geq -(m + n) + 1} \supset \cdots \supset W_{\geq -m + 1} \supset W_{\geq -m} \supset \emptyset.  
\end{equation}
Each graded piece of the filtration has dimension exactly $1$.  

Now, note that the trivial $H$-representation appears in the representation $(n,m)$ if and only if both $n$ and $m$ are even: the first condition guarantees that the zero weight appears in the representation $(n,m)$ and the second condition insures that $\pm Id \subset H$ acts trivially.  Then, by Lemma \ref{lem:mult} combined with Theorem \ref{thm:affinefp}, the multiplicity of the representation $(n,m)$ in $k[Sym^2 V^{\vee}]$ is identified with the space of linear maps $f$ from $W$ to the $1$-dimensional $k$-vector space $k$, which preserve the filtrations and which are $H$-module homomorphisms.  

We can restrict ourselves to the case where both $n$ and $m$ are even.  Here the trivial representation of $H$ appears with multiplicity one.  By Schur's lemma, the space of $H$-module maps from $(n,m)$ to the trivial representation is exactly $1$-dimensional.  In a basis of the vector space underlying $(n,m)$ in which $H$ acts diagonally, this family of linear maps is given parametrically by the vector $(0,\ldots,0,\lambda,0,\ldots 0)$; this is a vector of length $2n+1$ and $\lambda$ occurs at the $n$-th position.  Now, change bases to a basis in which the filtrations are given by coordinate subspaces; for the representation $(n,m)$, this involves conjugating by the $n$-th symmetric power of the matrix $\begin{pmatrix} 1 & 1 \\ 0 & 1\end{pmatrix}$.  Asking that the induced linear map preserves the filtrations constrains $m$ to be positive.  

Combining the two conditions discussed above, the representation $(n,m)$ appears in the coordinate ring $k[Sym^2 V]$ if and only if $n$ is even and $m$ is positive and even; when it appears, it appears with multiplicity $1$.
\end{ex}

\begin{ex}
\label{ex:matrices}
Finally, consider $GL_2 \times GL_2$ acting on the space $M_2(k)$ of $2 \times 2$-matrices over $k$ by $(g_1,g_2) \cdot A = g_1 A g_2^{-1}$.  We now completely describe the $GL_2 \times GL_2$-module structure of $k[M_2(k)]$.  The $2 \times 2$ identity matrix $Id \in M_n(k)$ gives a $k$-point in the open $GL_2 \times GL_2$-orbit.  The stabilizer of $Id$ is the diagonal subgroup $GL_2 \subset GL_2 \times GL_2$ and we can identify the open orbit with $(GL_2 \times GL_2)/GL_2$; the complement of the open orbit is the set of matrices $A$ satisfying $det A = 0$.  The Bruhat decomposition shows that $GL_2 \times GL_2/GL_2$ is spherical and hence that $M_2(k)$ is spherical.  The space $M_2(k)$ has three $GL_2 \times GL_2$-orbits; the open orbit, a codimension $1$ orbit corresponding to matrices with rank exactly $1$, and the zero matrix.  

We set
\begin{equation}
\mu(t) = \begin{pmatrix} t & 0 \\ 0 & t \end{pmatrix} \times \begin{pmatrix} 1 & 0 \\ 0 & t^{-1} \end{pmatrix}.
\end{equation}
Then note that $\lim_{t \rightarrow 0} \mu(t) \cdot Id$ exists and is a matrix of rank exactly $1$.  This one parameter subgroup satisfies the conclusion of Proposition \ref{prop:onepsg}

Every finite dimensional irreducible representation $W$ of $GL_2 \times GL_2$ is of the form $V_1 \boxtimes V_2$ where $V_1$ and $V_2$ are finite dimensional irreducible representations of the left and right hand factors of $GL_2$ respectively.  We use the characterization of finite dimensional representations of $GL_2$ given in Example $\ref{ex:binquadform}$.  Let $diag(t_1,t_2) \times diag(t_3,t_4)$ be a maximal torus in $GL_2 \times GL_2$.  We can compute the weights of an irreducible $GL_2 \times GL_2$-representation using the description of the weights of an irreducible $GL_2$-representation given in Equation \ref{eqn:weights}.  

We can compute the filtration on the vector space $W$ underlying a pair of representations indexed by the pairs of integers $(n,m)$ and $(n',m')$.  We have (see Equation \ref{eqn:constantfiltration}):
\begin{equation}
\label{eqn:filtration2}
\begin{matrix}
W \supset&  W_{\geq - (n + 2m - m')} \supset & W_{\geq -(n + 2m - (m' + 1)) } \supset & \cdots  \\
\supset W_{\geq -(n + 2m - (m' + n' - 1))} \supset &  W_{\geq -(n + 2m - (m' + n'))} \supset & \emptyset & \\
\end{matrix}.
\end{equation}
Each graded piece associated to this filtration has dimension exactly $n+1$.

Now, we use Lemma \ref{lem:mult} together with Theorem \ref{thm:affinefp} to identify the multiplicity of $W$ in $k[M_2(k)]$ with the space of $k$-linear maps $W \longrightarrow k$ which are $GL_2$-module homomorphisms for the diagonal copy of $GL_2$ and which preserve the filtration described in \ref{eqn:filtration2}.  A map can preserve the filtration if and only if $n+2m - (m' + n') \geq 0$.  By the Clebsch-Gordon formula, the restriction of $W$ to the diagonal copy of $GL_2$ contains the trivial representation if and only if $n = n'$ and $m = m'$; moreover, if it contains the trivial representation it contains it with multiplicity one.  

In summary, the representation of $GL_2 \times GL_2$ with weights $(n,m), (n',m')$ appears in $k[M_2(k)]$ if and only if $n = n'$, $m = m'$ and $m \geq 0$ in which case it appears with multiplicity one.  
\end{ex}

\begin{footnotesize}
\bibliographystyle{alpha}
\bibliography{asok-equivvb}
\end{footnotesize}

\end{document}